\documentclass[12pt]{amsart}
\newfont{\sheaf}{eusm10 scaled\magstep1}
\usepackage{amssymb}
\usepackage{showkeys}

\newcommand{\Proof}{{\it Proof. }}
\newcommand{\QED}{{\hfill $Q.E.D.$}}

\newtheorem{teo}{Theorem}[section]
\newtheorem{df}[teo]{Definition}
\newtheorem{lem}[teo]{Lemma}

\newtheorem{oss}[teo]{Remark}

\newtheorem{conj}[teo]{Conjecture}
\newcommand\sB{{\mathcal B}}

\newcommand\sK{{\mathcal K}}

\newcommand\De{\Delta}

\def\Bbb{\bf}

\def\C{{\Bbb C}}

\newcommand{\CC}{\ensuremath{\mathbb{C}}}
\newcommand{\RR}{\ensuremath{\mathbb{R}}}
\newcommand{\ZZ}{\ensuremath{\mathbb{Z}}}
\newcommand{\QQ}{\ensuremath{\mathbb{Q}}}

\newcommand{\MMM}{\ensuremath{\mathcal{M}}}

\newcommand{\NNN}{\ensuremath{\mathcal{N}}}
\newcommand{\NN}{\ensuremath{\mathbb{N}}}
\newcommand{\hol}{\ensuremath{\mathcal{O}}}

\newcommand{\PP}{\ensuremath{\mathbb{P}}}

\newcommand{\ra}{\ensuremath{\rightarrow}}

\def\eea{\end{eqnarray*}}
\def\bea{\begin{eqnarray*}}

\newcommand\dual{\mathrel{\raise3pt\hbox{$\underline{\mathrm{\thinspace d
\thinspace}}$}}}

\newcommand\qe{\ifhmode\unskip\nobreak\fi\quad $\Box$}       % box for QED

\def\BOX{\hfill\lower.5\baselineskip\hbox{$\Box$}}
\input xy
\xyoption{all}
\CompileMatrices

\title[ALGEBRAIC SURFACES: MODULI, REAL, DIFF., SYMPL. STRUCTURES]
{ALGEBRAIC SURFACES AND THEIR MODULI SPACES: REAL,
   DIFFERENTIABLE AND SYMPLECTIC STRUCTURES}
\author{Fabrizio Catanese, Universit\"at Bayreuth}

\begin{document}
\thispagestyle{empty}
%\LARGE
%\baselineskip=1cm

\maketitle

\section{RIPOSTE ARMONIE?}

At the onset of  the theory of algebraic surfaces (ca. 1870-1895)
Noether, Enriques and Castelnuovo said that while algebraic curves
had been  made by God, algebraic surfaces were made by the devil.

But Enriques, in the final page 464 of his famous book
` Le Superficie Algebriche' ( \cite{enriques}, published posthumously in
1949), says that really God   created  for algebraic surfaces a
higher level order
of `hidden harmonies' , where  an incredible beauty  shines.  The richness
and beauty of their properties, which were only after a long time and hard work
revealed, should inspire not only a sense of awe in the contemplation
of this divine order, but also the hope for researchers that
difficulties, doubts and
contradictions which pave their way to discovery are eventually to fade away
and uncover this divine light of harmony.

What does it mean that everything goes `right' for curves?
In order to explain this, we must first of all  recall the
   basic definition of a complex projective variety.

But why 'in primis' do we need  complex projective varieties?

The point is that in general, in `real life' (i.e., in applied
science), one wants to solve
polynomial equations with real coefficients, and find real
solutions. But since by the fundamental theorem of algebra every
`univariate' polynomial $P(x) \in \CC [x]$ (a polynomial in a single 
variable $x$)
of degree $n$ has exactly $n$ complex roots (counted with multiplicity),
the simple but basic approach is to view the real numbers as the subset of
the complex numbers fixed by complex conjugation,
hence the first approach is to start to look at
the complex solutions, and only later
to look at the action of complex conjugation
on the set of complex solutions.

Given moreover a system of polynomial equations in several variables
$( z_1, ... z_n)$, in order to have
some continuity of the dependence of the solutions upon the choice
of the coefficients, one reduces it to a system of
homogeneous equations,
$$ f_i(z_0, z_1, ... z_n) = 0, i= 1, \dots r. \\$$
defining an algebraic set or  a projective Variety
$$X  \subset  \PP^n_\CC,  X : = \{ z : = (z_0, z_1, ... z_n)| f_i(z)
= 0 \ \forall i = 1, \dots r \}$$
(one finds then the solutions  originally sought for by
setting $ z_0 = 1$).

We assume throughout here that $X$ is a smooth complex projective
variety, i.e., that  $X$ is a smooth
compact submanifold of $\PP^n : =  \PP^n_{\CC}$  of complex dimension $d$,
and that $X$ is connected.

Observe that $X$ is a compact oriented real manifold of real
dimension $2d$, so, for instance, a complex surface gives
rise to a real 4-manifold.

For $d=1$, $X$ is a complex  algebraic curve (a real surface, called
also Riemann surface) and its basic invariant
is the genus $g= g(X)$.

The genus $g$ is defined as
the dimension of the vector space $H^0 (\Omega^1_X)$ of rational 1-forms
$$\eta = \sum_i \phi_i(z) dz_i $$  which are homogeneous of degree
zero and are regular, i.e.,  do
not have poles  on $X$ (the  $\phi_i(z)$'s are here rational functions of $z$).

It turns out that the genus determines the topological and the differentiable
manifold underlying $X$: its intuitive meaning is the `number of handles'
that one must add to a sphere in order to obtain $X$ as a topological space.

Actually, as conjectured by Mordell and proven by Faltings
(\cite{faltings1},\cite{faltings2}, see also \cite{bom-mord} for
an `elementary' proof),
it also governs the arithmetic aspects of $X$: if  the coefficients of the
polynomial equations defining $X$ belong to  $\QQ$,
or more generally to a number field $k$ (a finite extension of $\QQ$),
then  the number of solutions with coordinates  in $\QQ$ (respectively:
with coordinates in $k$) is finite if
the genus $ g = g (X)  \geq 2$.

The rough classification of curves is the following :

\begin{itemize}

\item $g=0$ : $X \cong \PP^1_{\CC}$, topologically $X$ is a sphere $S^2$
of  real dimension $2$.
\item $g=1$ : $X \cong \CC / \Gamma$ , with $\Gamma$ a discrete subgroup
$\cong \ZZ^2$: $X$ is called an elliptic curve, and topologically we have
a real 2-torus $ S^1 \times S^1$.
\item $ g\geq 2 $ : then we have a `curve of general type', and topologically
we have a sphere with $g$ handles.
\end{itemize}

Moreover, $X$ admits a metric of constant curvature, positive if $g=0$,
negative if $ g\geq 2$, zero if $ g=1$.

Finally, we have a {\bf Moduli space} $\frak M_g$, an open set of a
complex projective
variety, which parametrizes the isomorphism classes of compact complex curves
of genus $g$. $\frak M_g$ is connected, and it has  complex dimension
$ (3g-3)$ for $ g \geq 2$.

Things do not seem so devilish when one learns that also for algebraic surfaces
of general type there exist similar moduli spaces $\frak M_{x,y}$
(by the results of \cite{bom}, \cite{gieseker}).

Here, $\frak M_{x,y}$ parametrizes isomorphism classes of minimal (smooth
projective) surfaces
of general type $S$ such that $ \chi (S) = x, K^2_S = y$.

Again, these two numbers are determined algebraically, through the
dimensions of
certain vector spaces of differential forms without poles,
namely, we have $ \chi (S) : = 1 - q(S) + p_g(S) , K^2_S : = P_2(S)- \chi(S)$,
where:
$$q(S): = dim_{\CC} H^0 (\Omega^1_S), p_g(S): = dim_{\CC}H^0 (\Omega^2_S),
P_2(S): = dim_{\CC}H^0
({\Omega^2_S}^{\otimes 2}).$$

As does the genus of an algebraic curves, these  numbers are determined
by the topological structure of $S$.

{\bf BAD NEWS : WE SHALL SEE THAT THESE TWO NUMBERS DO
NOT DETERMINE THE TOPOLOGY  of $S$!}

{\bf GOOD NEWS : $\frak M_{x,y}$ HAS FINITELY MANY CONNECTED COMPONENTS!}

The above finiteness statement is good news because
the connected components of $\frak
M_{x,y}$ parametrize

{\em Deformation classes
of surfaces of general type}, and, by a classical theorem of Ehresmann
(\cite{ehre}),
deformation equivalent varieties are diffeomorphic.

Hence, fixed the two numerical invariants $\chi (S) = x, K^2_S = y$,
which are determined
by the topology of $S$ (indeed, by the  Betti numbers of $S$),
we have a finite number of
differentiable types.

In the next section we shall  accept the bad news, trying to learn something
from them.

\section{Field automorphisms, the absolute Galois group and conjugate
varieties.}

Let $\phi : \CC \ra \CC$ be a field automorphism:
then, since $$\phi (x) = \phi ( 1 \cdot x) =\phi (1) \cdot \phi (x)
,$$ it follows
that $\phi (1) = 1 $, therefore $\phi (n) = n  \ \forall n \in \NN$,
and also $\phi |_{\QQ}  = Id_{\QQ}$.

We recall the first algebra exercises, quite  surprising for  first
semester students:
the real numbers have no other automorphism than the identity,
while the complex numbers have
`too many'.

\begin{lem}
$ Aut (\RR) = \{ Id \}$
\end{lem}

\Proof
For each choice of $ x , a \in \RR$, $\phi (x + a^2) = \phi (x) +
\phi (a)^2$, thus $\phi$
of a square is a square,
$\phi$ carries the set of squares $\RR_+$ to itself,  $\phi $ is increasing.
But $\phi$ equals the identity on $\QQ$: thus $\phi$ is the identity.
\QED

On the other hand,  the theory of transcendence bases
and the theorem of Steiniz  (any bijection between
two transcendence bases $\sB_1$ and $\sB_2$ is realized
a suitable automorphism) tell us:
\begin{lem}
$ | Aut (\CC) | = 2^{2^{\aleph_0}}$
\end{lem}

\begin{oss}
Paolo Maroscia informed me, after the lecture,
that the above result is the content
of a short note of Beniamino Segre (\cite{segre}) published 60 years ago.

Observe now that the only continuous automorphisms of $\CC$  are
the identity and the complex conjugation
$\sigma$, such that $\sigma (z) := \bar{z} = x - i y$. All the others
are impossible (just very hard?) to visualize!

2) The field of algebraic numbers $\overline{\QQ}$ is the
subfield of $\CC$,
$\overline{\QQ} : = \{ z \in \CC | \exists P \in \QQ[x]
s.t. \ P (z) = 0 \}.$   It is carried to itself by any field automorphism
of $\CC$.

The fact that $ Aut (\CC)$ is so large is essentially due to the fact that
the kernel of $Aut (\CC) \ra Aut (\overline{\QQ})$ is very large.

The group $Aut (\overline{\QQ})$ is called the {\bf absolute Galois group}
   and denoted
by
$ Gal (\overline{\QQ}, \QQ)$. It is one of the most interesting objects of
investigation in algebra and arithmetic.

Even if we have a presentation %(\cite{fried})
of this
group, still our information about it is quite scarce. A presentation of
a group $G$ often does not even answer the question:
is the group $G$ nontrivial? The solution to this question is often
gotten if we
   have a {\bf representation} of the group $G$, for instance, an action
of $G$ on a set $M$ that
can be very well described.
\end{oss}

For instance, $M$ could be here a moduli space.

To explain how  $ Gal (\overline{\QQ}, \QQ)$ acts on several moduli spaces $M$
we come now to the
crucial notion of a conjugate variety.

\begin{oss}
1) $\phi \in Aut (\CC)$ acts on $ \CC [z_0, \dots z_n]$,
by sending $P (z) = \sum_{i =0}^n \ a_i z ^i
\mapsto  \phi (P) (z) : = \sum_{i =0}^n \ \phi (a_i) z ^i$.

2) Let $X$ be as above a projective variety
$$X  \subset  \PP^n_\CC,  X : = \{ z | f_i(z) = 0 \forall i \}.$$

The action of $\phi$ extends coordinatewise to $ \PP^n_\CC$,
and carries $X$ to another variety, denoted $X^{\phi}$,
and called the {\bf conjugate variety}. Since $f_i(z) = 0 $ implies
$\phi (f_i)(\phi (z) )= 0 $, we see that

$$  X^{\phi}  = \{ w | \phi (f_i)(w) = 0 \forall i \}.$$

\end{oss}

If $\phi$ is complex conjugation, then it is clear that the variety
$X^{\phi}$ that we obtain is diffeomorphic  to $X$, but in general,
what happens when $\phi$ is not continuous ?

For curves, since in general the dimensions of spaces of
differential forms of a fixed degree and without poles are the same
for $X^{\phi}$ and $X$, we shall obtain a curve of the same genus,
hence $X^{\phi}$ and $X$ are diffeomorphic.

But for higher dimensional varieties this breaks down,
as discovered by Jean  Pierre Serre in the 60's (\cite{serre}),
who proved the existence of a field automorphism  $\phi \in
Gal(\bar{\QQ} /\QQ)
$, and a variety $X$ defined over $\bar{\QQ}$ such that
$X$ and the Galois conjugate variety $X^{\phi}$  have
   non isomorphic fundamental groups.

In  work  in progress in collaboration with Ingrid Bauer and Fritz
Grunewald (\cite{bcgGalois}) we discovered  wide classes of algebraic surfaces
for which the same phenomenon holds.

Shortly said, our method should more generally yield
   a way to transform the bad
news into good news:

\begin{conj}
Assume that $\phi \in Gal(\bar{\QQ} /\QQ) $
is different from the identity and from complex
conjugation. Then there is a minimal surface of general type $S$ such that
$S$ and $S^{\phi}$ have non isomorphic fundamental groups.
In particular, $S$ and $S^{\phi}$ are not homeomorphic.

Hence the absolute Galois group $Gal(\bar{\QQ} /\QQ)$
acts faithfully
on the  set of connected
components of the (coarse) moduli space of minimal surfaces of general type,
$$ \frak M : = \cup_{x,y \geq 1} \frak M_{x,y}. $$
\end{conj}

Concerning the above two statements, we should observe
   that, while the absolute Galois group
$Gal(\bar{\QQ} /\QQ)$ acts on the set  of connected components of
$ \frak M $, it does not  act on the set of isomorphism
classes of fundamental groups of
surfaces of general type: this means that, given two varieties $X,
X'$ with isomorphic
fundamental groups, their conjugate varieties $X^{\phi}, {X'}^{\phi}$
do not need to have isomorphic fundamental groups.  Else,
not only  complex conjugation would not change the isomorphism class of the
fundamental group, but also the minimal normal subgroup containing it
(which is very
large) would do the same.

Let me end this section giving a few hints about the main ideas
and methods for our proposed approach, which depends on a single
general conjecture
about faithfulness of the action of the absolute Galois group
$Gal(\bar{\QQ} /\QQ)$
on the isomorphism classes of unmarked triangle curves.

An elementary but key lemma describes our candidate triangle curves
for the above conjecture.

   Fix a positive integer $g \in \NN$, $g \geq 3$, and define,
   for any complex number $ a \in \CC \setminus \{ -2g,0,1, \ldots , 2g-1 \} $,
$C_a$  as the hyperelliptic
curve   of genus $g$  given by the equation
$$ w^2 = (z-a) (z + 2g) \Pi_{i=0}^{2g-1} (z-i) .$$
We have then:
\begin{lem} Consider two complex numbers  $a,b$ such that $a \in  \CC
   \setminus \QQ $: then $C_a
\cong C_b$ if and only if $a = b$.
\end{lem}

Through the above lemma algebraic numbers are therefore encoded into
isomorphism
class of curves.

We use then the method of proof of the well known Belyi theorem:

\begin{teo} {\bf (BELYI)}
An algebraic curve $C$ can be defined over $\overline{\QQ}$  if and only
if there exists a holomorphic map $ f : C \ra \PP^1_{\CC}$ with branch  set
(set of critical values) equal to
  $\{0, 1, \infty \}$.
\end{teo}

Assume now that $a$ is algebraic, i.e., that $ a \in \overline{\QQ}$:
   take a Belyi function for
$C_a$ (i.e., $ f_a : C_a \ra \PP^1_{\CC}$ with branch
set  $\{0, 1, \infty \}$)
and its normal closure
$D_a \ra \PP^1_{\CC}$.
We have then constructed a triangle curve $D_a$ according to the following

\begin{df}
$D$ is a {\bf TRIANGLE CURVE} if there is a finite group $G$ acting effectively
on $D$ and with the property that $ D / G \cong \PP^1_{\CC}$, and the
quotient map
$ f : D \ra \PP^1_{\CC} \cong D / G $ has $\{0, 1, \infty \}$ as branch set.

A {\bf  MARKED TRIANGLE CURVE} is a triple $(D,G,i)$ where $D,G$ are as above,
and where  we have fixed an embedding $ i : G \ra Aut (D)$.

Two marked triangle curves $(D,G,i)$, $(D',G',i)$ are isomorphic iff there
exists isomorphisms $ D \cong D'$, $G \cong G'$ which transform
$ i $ into $ i'$.
\end{df}

Let us explain now the basic idea which lies behind our new results: 
the theory of
surfaces  isogenous  to a product, introduced in \cite{isogenous}(see also
\cite{cat4}), and which  holds more generally also for higher dimensional
varieties.

\begin{df}

1) A surface {\bf isogenous to a (higher) product}
is a compact complex surface $S$ which
  is  a quotient
$S = (C_1 \times C_2) /G$ of a product of curves of resp. genera
$\geq 2$ by the
free action of a finite group $G$.

2)A {\bf Beauville \ surface}   is a surface  isogenous to a (higher) product
which is {\bf rigid}, i.e., it has no nontrivial deformation. This
amounts to the condition,
in the case where the two curves $C_1$ and $C_2$ are nonisomorphic,
that  $ (C_i, G)$ is a triangle curve.

\end{df}

For surfaces isogenous to a product holds the following (\cite{isogenous},
\cite{cat4}):

\begin{teo}
Let $S = (C_1 \times C_2) /G$ be a surface isogenous to a product.
Then any surface $X$ with the
same topological Euler number and the same fundamental group as $S$
is diffeomorphic to $S$. The corresponding subset of the moduli space
$\frak M^{top}_S = \frak M^{diff}_S$, corresponding to surfaces homeomorphic,
resp, diffeomorphic to $S$, is either irreducible and connected
or it contains
two connected components which are exchanged by complex
conjugation.
\end{teo}

If $S$ is a Beauville surface this implies: $X \cong S$ or $X \cong \bar{S}$.
It follows also that a Beauville surface is defined over $\overline{\QQ}$,
whence the Galois group $Gal (\overline{\QQ} , \QQ)$ operates on the
discrete subset of the moduli space $\mathfrak M$ corresponding
to Beauville surfaces.

Work in progress with the same coauthors (Ingrid Bauer and Fritz Grunewald)
aims at proving also the following
\begin{conj}
The absolute Galois group $Gal (\overline{\QQ} , \QQ)$ operates
faithfully on the
discrete subset of the moduli space $\mathfrak M$ corresponding
to Beauville surfaces.

\end{conj}
Already established is the following (\cite{bcgGalois})
\begin{teo}
Beauville surfaces  yield explicit examples of conjugate surfaces with
nonisomorphic fundamental groups whose completions are isomorphic
(the completion of a group $G$ is the inverse limit
$$ {\hat G} = lim_{K\subset G \ normal \ of \ finite \ index}  \ \ 
(G/K)  \ ).$$
\end{teo}
Our candidate triangle curves  $C_a$ determine now a family $\frak N_a$
consisting of
all the possible surfaces isogenous to a product of the form
$S : = (D_a \times D') / G$, where the genus of $D'$ is fixed, and $G$
acts without fixed points on $D'$.

Using the theory of surfaces {\bf isogenous to a product}, referred to above,
follows easily that:

1) $\frak N_a$ is a union of connected components of $\frak M$

2) $ \phi (\frak N_a) = \frak N_{\phi (a)}$.

Assume that for each $\phi \in Aut (\overline{\QQ})$ which is nontrivial
we can find $ a$ such that, setting $ b : = \phi (a) $:

3) $ a \neq b$ and $\frak N_a$ and $\frak N_b$ do not intersect.

The desired conclusion would then be  that,  since $ \phi (\frak N_a) $
and $\frak N_a$ do not intersect by 2),3), hence
$\phi$ acts nontrivially on the set of connected components
of $\frak M$.

The condition that $\frak N_a$ and $\frak N_b$ intersect easily implies,
by the structure theorem for surfaces isogenous to a product, that
the two triangle curves $D_a$ and $D_b$ are isomorphic.
There is thus an isomorphism  $ F : D_a \ra D_b$ which transforms the
action of $G_a$ on $D_a$ into the action of $G_b$ on $D_b$.
Identifying  $G_a$ with $G_b$ under the transformation $\phi$, one sees
however that $F$ is only `twisted' equivariant.  This means that there
is an isomorphism $\psi \in Aut (G)$ such that $ F ( g (x)) = \psi (g) (x)$.

  If $(D_a, G, i_a)$ and $(D_b, G, i_b)$ are isomorphic as marked triangle
curves (for instance, if $\psi$ is inner), then it follows that $C_a$ 
is isomorphic
to $C_b$,  and we derive   a contradiction, that $ a=b$.
In other words, our previous lemma shows that the
absolute Galois group $Gal(\bar{\QQ} /\QQ)$
acts faithfully
on the  set of isomorphism classes of marked triangle curves.

The main point is to find
such an $a \neq b = \phi (a)$ with the above  property that the
group $G$ has only inner automorphisms.

Indeed, the only crucial property which should be proven amounts to
the following

\begin{conj}
The absolute Galois group $Gal(\bar{\QQ} /\QQ)$
acts faithfully
on the  set of isomorphism classes of (unmarked) triangle curves.
\end{conj}

There are other  interesting open  problems:

{\bf Question 1.} Existence and classification of Beauville surfaces, i.e.,

a) which finite groups $G$ can occur?

b) classify all possible Beauville surfaces for a given finite group $G$.

We have made a substantial progress (\cite{bcg}, \cite{almeria})
on question a) which leads
for instance to substantial evidence in the direction of
the following

\begin{conj} Every finite nonabelian simple group occurs except
$\frak U_5$.
\end{conj}

\section {Surfaces of general type,  DEF= DIFF ? and beyond.}

Let $S$ be a minimal surface of general type: then we saw that to $S$
we attach two positive integers $\geq 1$
$$ x = \chi (\hol_S), \ y
=K^2_S$$

which are invariants of the oriented topological type
of $S$ (they are determined by the Euler number $e(S)$ and
by the signature $\tau (S)$).

The {\bf moduli space} $\frak M_{x,y}$ of the surfaces with
invariants $(x,y)$ is a
quasi-projective variety defined over the integers, and it has a finite number
of irreducible components.

For fixed $(x,y)$ we have thus a finite number of possible differentiable
types, and a fortiori a finite number of topological types.

Michael  Freedman's big Theorem of 1982 (\cite{free}) shows that there are
indeed at most  two topological structures if moreover the  surface $S$ is
assumed to be simply connected (i.e., with trivial fundamental group).

Topologically, our 4-manifold is then obtained  from  very simple
building blocks,
one of them being the K3 surface, where:

\begin{df}
A {\bf K3} surface is a smooth surface of degree 4 in $\PP^3_{\CC}$.
\end{df}

Observe moreover that a complex manifold carries a natural orientation
corresponding to the complex structure, and, in general, given an oriented
differentiable manifold $M$, $M^{opp}$ denotes the same manifold, but endowed
with the opposite orientation. This said, we can explain the corollary
of Freedman's theorem for the topological manifolds underlying
simply connected (compact) complex surfaces.

There are two cases which are distinguished as follows:

\begin{itemize}
\item
   $S$ is EVEN, i.e., its intersection form on $H^2 (S, \ZZ)$ is even:
then $S$ is a connected sum of copies of

   $ \PP ^1_{\CC} \times \PP ^1_{\CC}$ and of a $K3$ surface
if the signature $\tau(S)$ is negative, and of copies of
$ \PP ^1_{\CC} \times \PP ^1_{\CC}$ and of a $K3$ surface with
reversed orientation if the signature is positive.

\item
$S$ is ODD : then $S$ is a connected sum of copies of
$\PP ^2_{\CC} $ and ${\PP ^2_{\CC}}^{opp} .$
\end{itemize}

We recall that the connected sum is the operation which,
from two oriented manifolds $M_1$ and $M_2$,
glues together the complements of two open balls $ B_i \subset M_i$,
having a differentiable (resp. : tame) boundary, by identifying
together the two boundary spheres via an orientation
reversing diffeomorphism.

Kodaira and Spencer defined quite generally (\cite{k-s58}) for  compact complex
manifolds
$X, X'$  the equivalence relation called {\bf deformation equivalence}: this,
   for  surfaces of general type, means that the corresponding
isomorphism classes
   yield points in the same connected component
of the moduli space $\frak M$.
The cited theorem of Ehresmann guarantees that DEF $\Rightarrow $ DIFF:

\begin{oss}
Deformation equivalence implies the existence of a diffeomorphism carrying
the canonical class $K_X$ to the canonical class $K_{X'}$.
\end{oss}

In the 80' s , groundbreaking work of Simon Donaldson (\cite{don1},
\cite{don2},\cite{don3},\cite{don4}, see also \cite{d-k}) showed that
homeomorphism and diffeomorphism differ drastically
for projective surfaces.

\begin{oss}
A refinement of Donaldson's theory, made by Seiberg and Witten
(see \cite{witten},\cite{don5},\cite{mor}),
showed then more easily that  a diffeomorphism $\phi : S \ra S'$
between minimal surfaces of general type satisfies
$\phi^* (K_{S'}) = \pm  K_S$.
\end{oss}

Based on the successes of gauge theory,
the following conjecture was made (I had been writing five years before
the opposite conjecture, in \cite{katata}, but almost no one believed in it):

FRIEDMAN-MORGAN'S SPECULATION (1987) (see \cite{f-m1}):  DEF =  DIFF
(Differentiable equivalence and deformation equivalence coincide for surfaces)

However, finally the question was answered negatively in every
possible way (\cite{man4},\cite{k-k},\cite{cat4},\cite{c-w},\cite{bcg},
see also \cite{catcime} for a  rather comprehensive survey):

\begin{teo}
(Manetti '98, Kharlamov -Kulikov 2001, C. 2001, C. - Wajnryb 2004,
Bauer- C. - Grunewald 2005 )

The Friedman- Morgan speculation
does not hold true.
\end{teo}
\begin{itemize}
\item
(1) Manetti used $(\ZZ /2)^r$-covers of  blow ups of the quadric
$Q : =  \PP ^1_{\CC} \times \PP ^1_{\CC}$, his
surfaces have $b_1 =0$, but are not simply connected.
\item
(2) Kharlamov and Kulikov  used   quotients $S$ of the unit ball in $\CC^2$:
the surfaces they use are rigid but with infinite fundamental group.
\item
(3) I used non rigid surfaces isogenous to a product
   $S = ( C_1 \times C_2 )/ G$, thus with $b_1 > 0$ and a fortiori
the surfaces have infinite fundamental group.
\item
(4) The examples given with Bauer and Grunewald are Beauville
surfaces, again
the surfaces are rigid, thus they have $b_1 = 0$ but  infinite 
fundamental group.
\item
(5) The examples obtained with Wajnryb are
instead simply connected, i.e., they have trivial fundamental group.

\end{itemize}

\bigskip

{\bf Common feature of (2),  (3) and (4) : we take $S$ and the
conjugate surface $\bar{S}$
( thus  $\bar{S}
, S$ are diffeomorphic), and if $\bar{S}$
and $S$  were deformation equivalent, there
would be a self-diffeomorphism $\psi$ of $S$ with $\psi^* (K_{S}) = - K_S$.
    If $\psi$  exists,  it should  be antiholomorphic ( by general properties
of these surfaces).
The technical heart of the proof is to construct examples where this
cannot happen,
and the fundamental group is heavily used for this issue.}

\bigskip

After the first counterexamples were found, the following
weaker conjectures were
posed:

\section{ Weakenings of the conjecture by Friedman and Morgan }
\begin{itemize}
\item
(I) require a diffeomorphism $\phi : S \ra S'$ with
$\phi^* (K_{S'}) =   K_S$.
\item
(II) require the surfaces to be simply connected (1-connected).

\end{itemize}

Even these weaker conjectures were disproven in my joint work with
Bronek Wajnryb (\cite{c-w}), which I will now briefly describe.

The  simply connected examples we used, called {\bf abc-surfaces}, are
a special case of a class of surfaces which I introduced in 1982
(\cite{cat1}), namely, bidouble covers of the quadric
and their natural deformations.

Bidouble covers of the quadric are smooth projective complex surfaces $S$
endowed with a (finite) Galois
covering $ \pi : S \ra Q : = \PP^1 \times \PP^1$ with Galois group
$(\ZZ /2 \ZZ)^2$.

More concretely, they are defined by a single pair of equations

\begin{eqnarray*}
   & z^2= f_{(2a,2b)}(x_0 , x_1; y_0 , y_1)\\
   & w^2=g_{(2c,2d)}(x_0 , x_1; y_0 , y_1)\\
\end{eqnarray*}

where $ a,b,c,d\in\NN^{\geq3}$ and the notation $f_{(2a,2b)}$ denotes that
$f$ is a bihomogeneous polynomial, homogeneous of degree $2a$ in the
variables $x$,
and of degree $2b$ in the variables $y$.

These surfaces are simply connected and minimal of general type,
and they  were introduced in \cite{cat1} in order to show that the
   moduli spaces
$\frak M_{\chi,K^2}$ of smooth minimal surfaces of general type $S$ with
$K^2_S = K^2 , \chi(S): = \chi (\hol_S) = \chi$
need not be equidimensional or irreducible (and indeed the same holds
for the open and closed subsets $\frak M_{\chi,K^2}^{00}$ corresponding to
simply connected surfaces).

Given in fact our four integers $ a,b,c,d\in\NN^{\geq3}$,  considering
the socalled natural deformations of these bidouble covers,
defined by equations \footnote{in the following formula,
a polynomial of negative degree is identically zero. }
\begin{eqnarray*}
   & z^2= f_{(2a,2b)}(x_0 , x_1; y_0 , y_1) + w \
\Phi _{(2a-c,2b-d)}(x_0 , x_1; y_0 , y_1)\\
   & w^2=g_{(2c,2d)}(x_0 , x_1; y_0 , y_1) + z \
\Psi _{(2c-a,2d-b)}(x_0 , x_1; y_0 , y_1)\\
\end{eqnarray*}
one defines a bigger open subset $\NNN'_{a,b,c,d}$ of the moduli space,
whose closure $\overline{\NNN'}_{a,b,c,d}$
is an irreducible component of $\frak M_{\chi,K^2}$, where
   $ \chi = 1 + (a-1)(b-1) + (c-1)(d-1) + (a+c-1)(b+d-1)$, and $K^2 =
(a+c-2)(b+d-2)$.

The  {\bf abc-surfaces} are obtained as the special case where $ b = d$,
and the upshot is that, once the  values of the integers  $b$
and  $a+c$ are fixed, one obtains diffeomorphic
surfaces.

In other more technical words the  {\bf abc-surfaces}
   are the {\bf natural deformations} of
    $({\ZZ}/2)^2$-covers of
   $(\PP^1\times
\PP^1)$, of {\em simple type}  $(2a,2b),(2c,2b)$,
which means that they are defined by equations

$$  (***) \ \ \  \ z_{a,b}^2 =  f_{2a,2b} (x,y)  + w_{c,b} \phi_{2a-c,b}(x,y)$$
    $$ w_{c,b}^2 = g_{2c,2b}(x,y) +  z_{a,b}  \psi_{2c-a,b}(x,y)$$

where f,g ,$\phi, \psi$, are bihomogeneous polynomials ,  belonging to
     respective vector spaces of sections of line bundles:

     $ f \in H^0({\PP^1\times \PP^1}, {\hol}_{\PP^1\times \PP^1}(2a,2b)) $,

$ \phi \in H^0({\PP^1\times \PP^1}, {\hol}_{\PP^1\times \PP^1}(2a-c,b))$ and

$ g \in H^0({\PP^1\times \PP^1}, {\hol}_{\PP^1\times \PP^1}(2c,2d))$,

$\psi \in H^0({\PP^1\times \PP^1}, {\hol}_{\PP^1\times \PP^1}(2c-a,b))$.

\bigskip

The   main new result of \cite{c-w} is the following

\begin{teo}\label{diffabc} {\bf (C. -Wajnryb )}
Let $S$ be an $(a,b,c)$ - surface and $S'$ be an
$(a+1,b,c-1)$-surface. Moreover, assume that $a,b,c-1 \geq 2$. Then
$S$ and $S'$ are diffeomorphic.
\end{teo}

This result  couples then with a more technical result:

\begin{teo} \label{nondef}  {\bf (C. -Wajnryb )}
Let  $S$, $S'$ be simple bidouble
covers of  ${\PP}^1
\times  {\PP}^1  $ of respective
     types ((2a, 2b),(2c,2b), and (2a + 2k, 2b),(2c - 2k,2b) , and assume

\begin{itemize}
\item
     (I) $a,b,c, k$
are strictly  positive even integers with $ a, b, c-k \geq 4$
\item
(II)  $ a \geq 2c + 1$,
\item
(III) $ b \geq c+2$  and either
\item
(IV1) $ b \geq 2a + 2k -1$ or\\
(IV2) $a \geq b + 2$
\end{itemize}
Then  $S$ and  $S'$ are not deformation equivalent.
\end{teo}

    The  second theorem uses techniques which have been developed
in a long series of papers by the  author and by Marco Manetti in
a period of almost 20 years.

They use essentially

i) the local deformation
theory \'a la Kuranishi, but for the canonical models,

ii) normal
degenerations of smooth
surfaces and a study of
quotient singularities of rational double points and of their smoothings.

A detailed expositions for both theorems can be found
in the   Lecture Notes of the C.I.M.E. courses
`Algebraic surfaces and symplectic 4-manifolds' (see especially 
\cite{catcime}).

The result in differential topology obtained with Wajnryb is based instead on a
refinement of Lefschetz theory obtained by Kas (\cite{kas}).

This refinement allows us to encode the differential topology of a 4-manifold
$X$ Lefschetz fibred over $\PP^1_{\CC}$ (i.e.,
$f : X \ra PP^1_{\CC}$ has the property that all the fibres are smooth and
connected, except for a finite number which have a nodal singularity)
into an equivalence class of factorization of the identity
in  the Mapping class group $\MMM ap_g$ of a compact curve
$C$ of genus $g$.

The mapping class group, introduced by Max Dehn(\cite{dehn})
in the 30's, is defined for each manifold $M$
as
$$ \MMM ap (M) : =  Diff (M) / Diff^0 (M), $$
where $Diff^0 (M)$ is the connected component of the identity,
the so called subgroup of the diffeomorphisms which are
{\bf isotopic to the identity}.

A major advance in our knowledge of $\MMM ap_g$
was made by Hatcher and Thurston (\cite{h-t}), and the simplest known
presentation of this group is due to Wajnryb (\cite{wajnryb}, \cite{waj2}).

Verifying isotopy of diffeomorphisms is a difficult and very geometric task,
which is accomplished in \cite{c-w} by constructing chains of loops in the
complex curve $C$, which lead to a dissection of $C$ into open cells.
One has then to choose several associate Coxeter elements
to express a given diffeomorphism, used for glueing two
manifolds with boundary $M_1$ and $M_2$
in two different ways, as a product
of certain Dehn twists. This expression shows that this diffeomorphism extends
to the interior of $M_1$, hence that the results of the two glueing operations
yield diffeomorphic 4-manifolds.

I would like to finish  in the next section commenting
on some very interesting
open problems.

To discuss them, I need to explain
the connection with the theory of symplectic manifolds.

\section{Symplectic manifolds}

\begin{df}
A pair $(X,\omega)$ of a real manifold $X$, and of a real differential
2- form $\omega$ is called a {\bf Symplectic pair} if

i) $\omega$ is a symplectic form, i.e., $\omega$ is closed ( $ d \omega  = 0$)
    and  $\omega$ is nondegenerate at each point (thus $X$ has even
real dimension).

A symplectic pair $(X,\omega)$ is said to be {\bf integral} iff the De Rham
cohomology class of $\omega$ comes from $ H^2 ( X, \ZZ)$, or equivalently,
there is a complex line bundle $L$ on $X$ such that $\omega$ is a first Chern
form of $L$.

An {\bf almost complex structure} $J$ on $X$ is a differentiable endomorphism
of the real tangent bundle of $X$ satisfying $ J^2 = - 1$. It is said to be

ii) {\bf compatible with $\omega$ } if
$$ \omega (Jv, Jw) = \omega (v, w)$$

iii) {\bf tame }
if the quadratic form
$\omega ( v, J v)$ is strictly positive definite.

Finally, a symplectic manifold is a manifold admitting a symplectic form
$\omega$.

\end{df}

For long time (before the celebrated examples
of Kodaira and Thurston, \cite{kod}, \cite{th}) the basic examples of 
symplectic
manifolds were given by symplectic submanifolds
of the flat space $\CC^n$ and of K\"ahler manifolds,
in particular of the projective space $\PP^N_{\CC}$,
which possesses the Fubini-Study form

\begin{df}
The Fubini-Study form is the differential 2-form
$$\omega_{FS} := \frac{i}{2 \pi } \partial
\overline{\partial} log |z|^2,$$ where $z$ is the homogeneous coordinate
vector representing a point of $\PP^N $.

    In fact the above 2- form on $\CC^{N+1} \setminus \{0\}$ is invariant

1) for the action of $U(N, \CC)$ on homogeneous coordinate vectors

2) for multiplication of the vector $z$ by a nonzero holomorphic
scalar function $f(z)$ ($z$ and $ f(z) z$ represent the same point in
$\PP^N $), hence

3) $\omega_{FS}$ descends to a differential form on $\PP^N $,
being $\CC^*$-invariant.
\end{df}

\bigskip

I recently observed (\cite{cat02}, \cite{cat06}):

\begin{teo}
   A minimal surface of general type $S$ has a symplectic structure 
$(S,\omega)$,
unique up to symplectomorphism, and invariant for smooth deformation,
with class$(\omega) = K_S.$ This symplectic structure is called the
{\bf  canonical symplectic structure}.
A similar result holds for higher dimensional complex projective manifolds
with ample canonical divisor $K_X$.
\end{teo}

The above result is, in the case where $K_S$ is ample, a
rather direct consequence of the famous Moser's
lemma (\cite{moser}).
\begin{lem}
Let $f :\Xi \ra T$ a proper submersion of differentiable manifolds,
with $T$ connected,
and let $(\omega)$ be a 2-form on $\Xi$ whose fibre restriction
$\omega_t := \omega|_{X_t}$ makes each $X_t$ a
symplectic manifold.

If the class of $(\omega_t)$ in De Rham cohomology is constant, then the
$(X_t,\omega_t)$'s are all symplectomorphic.
\end{lem}

{\bf In the above theorem, when $K_S$ is ample,  it suffices to pull-back $1/m$
of the Fubini-Study metric by an
$m$-canonical embedding.}

In the case where $K_S$ is not ample, the proof
is more involved and uses techniques from
the following symplectic approximation theorem

\begin{teo}\label{glueing}
Let $X_0 \subset \PP^N$ be a projective variety with isolated
    singularities each admitting a smoothing.

Assume that for each singular point $x_h \in X$,
we choose a smoothing component $T_{j(h)}$ in the basis of the
semiuniversal deformation of the germ $(X, x_h)$.

Then (obtaining different results for each such choice) $X$ can be
approximated by symplectic submanifolds
$W_t$ of $ \PP^N$, which are diffeomorphic to the glueing
of the `exterior' of $X_0$ (the complement to the union $B = \cup_h B_{h}$ of
     suitable (Milnor) balls around the singular points) with
the  Milnor fibres $\MMM_h$, glued along the singularity links $\sK_{h,0}$.
\end{teo}
An important consequence is the following theorem (\cite{cat02}, \cite{cat06},
see also \cite{catcime} for a survey including the basics concerning the
construction of the Manetti surfaces)

   \begin{teo}
Manetti's surfaces  yield examples of surfaces of general type which are
not deformation equivalent but are canonically symplectomorphic.
\end{teo}

{\bf Questions:} 1) Are there (minimal) surfaces of general type which
are orientedly
    diffeomorphic through a diffeomorphism carrying the canonical class
to the canonical
class, but, endowed with their canonical symplectic structure,
are not
canonically   symplectomorphic?

2) Are there such   simply connected examples ?

3) Are the diffeomorphic abc-surfaces canonically symplectomorphic
(thus yielding a  counterexample to Can. Sympl = Def
also in the simply connected case)?

I am currently working with Wajnryb and L\"onne on the very difficult problem
of understanding the canonical symplectic structures of abc-surfaces 
(\cite{clw}).

To explain our result, let us go back to our equations

$$  (***) \ \ \  \ z_{a,b}^2 =  f_{2a,2b} (x,y)  + w_{c,b} \phi_{2a-c,b}(x,y)$$
    $$ w_{c,b}^2 = g_{2c,2b}(x,y) +  z_{a,b}  \psi_{2c-a,b}(x,y)$$

where $f,g$ are bihomogeneous polynomials as before, and instead we
allow $\phi,
\psi$, in the case where for instance the degree relative to $x$ is negative,
to be antiholomorphic.  In other words, we allow $\phi,
\psi$ to be sections of certain line bundles which are dianalytic
(holomorphic or antiholomorphic) in each variable $x, y$.

In this way we obtain a symplectic 4-manifold which (we call a
dianalytic perturbation and) is
canonically symplectomorphic to the bidouble cover we started with.
But now  we have gained that, for general choice of $f,g, \phi,
\psi$, the projection onto $\PP^1
\times
\PP^1$ is {\bf generic} and its branch curve $B$ (the locus of the 
critical values)
is a dianalytic curve with nodes and cusps as only singularities.

The only price one has to pay is to allow also negative nodes, i.e.,
nodes which in local holomorphic coordinates are defined by the equation
$$ ( y -\bar x) ( y + \bar x) = 0. $$

Now, projection onto the first factor $ \PP^1$ gives a movement of $n$
points in a fibre $ \PP^1$, which is encoded in the so called
{\bf vertical braid monodromy factorization}.

The first result that we have achieved  is the computation of this
vertical braid monodromy
factorization of the  above branch curve $ B \subset \PP^1 \times \PP^1$.

  The second very interesting result that we have obtained, and which is too
complicated to explain here in detail, is that certain invariants of 
these vertical
braid monodromy factorizations allow to reconstruct all the three numbers
$a,b,c$ and not only the numbers $b, a+c$, which determine the diffeomorphism
type.

This  result represent the first
positive step towards the realization
of a more general program set up by Moishezon (\cite{mois1}, \cite{mois2})
in order to produce braid monodromy invariants which should distinguish
connected components of the moduli spaces $\frak M _{\chi,K^2}$.

Moishezon's program is based on the consideration (assume here for simplicity
that $K_S$ is ample) of
a general projection $\psi_m : S \ra \PP^2$ of a pluricanonical
embedding $\psi_m : S \ra \PP^{P_m -1}$, and of the braid monodromy
factorization
corresponding to the (cuspidal) branch curve $\De_m$ of $\psi_m$.

An invariant of the connected component is here given by the equivalence class
(for Hurwitz equivalence plus simultaneous conjugation) of this
braid monodromy factorization. Moishezon, and later Moishezon-Teicher
calculated a coarser invariant, namely the fundamental group
$\pi_1(\PP^2- B_m)$. This group turned out to be not overly complicated,
and in fact, as shown in many cases in \cite{a-d-k-y}, it tends  to give
no extra information beyond the one given by the topological
   invariants of $S$ (such as $\chi, K^2$).

Auroux and Katzarkov showed instead (\cite{a-k}) that, for $ m >>0$, a
more general equivalence class (called
m-equivalence class, and allowing creation of a pair of neighbouring nodes,
one positive and one negative) of
the above braid monodromy factorization determines the canonical
symplectomorphism class of
$S$, and conversely.

The work by Auroux, Katzarkov adapted Donaldson's techniques for proving
the existence of symplectic Lefschetz fibrations (\cite{don6},\cite{donsympl})
in order to show  that each symplectic $4$-manifold is in a
natural way 'asymptotically' realized by  a generic symplectic covering of
$\PP^2_{\CC}$, given by almost holomorphic sections of a high multiple
$L^{\otimes m}$ of a complex line bundle $L$
whose class is the one of the given
integral symplectic form.

The methods of  Donaldson on one side, Auroux and Katzarkov on the other,
  use algebro geometric methods in order to produce invariants of
symplectic 4-manifolds.

For instance, in the case  of a generic symplectic covering of the
plane, we get a corresponding branch curve $\De_m$  which is a symplectic
submanifold with singularities only nodes and cusps.

To  $\De_m$ corresponds then a factorization in the braid group,
called m-th braid monodromy factorization: it  contains
only factors which are conjugates of $ \sigma_1^j$, not only with $ j
=  1, 2,3$
as in the complex algebraic case,
but also with $ j = -2$
(here $ \sigma_1$ is a standard half twist on a segment connecting two roots,
the first of the standard
Artin generators of the braid group).

Although the factorization is not unique (because it may happen that a pair of
two consecutive nodes, one positive and one negative, may  be created,
or disappear) one considers its m-equivalence class, and
   the authors show that this class, for $ m >> 1$, is an invariant
of the integral symplectic manifold.

In the case of abc-surfaces, consider now again the quadric  $ Q : = \PP^1
\times \PP^1$, and let $ p : S \ra \PP^2$ be the morphism
  obtained as the composition of
$\pi : S \ra Q$ with the standard (Segre) embedding
$ Q \hookrightarrow \PP^3$
and with a general projection  $\PP^3 \dashrightarrow \PP^2$.

In the special case of those particular abc-surfaces such that $ a+c = 2b$,
   the m-th pluricanonical mapping $\psi_m : S \ra \PP^{P_m -1}$
has a (non generic) projection given by the composition of $p$ with
a Fermat type map  $\nu_r : \PP^2  \ra \PP^2$
(given by $ \nu_r (x_0, x_1, x_2) = (x_0^r, x_1^r, x_2^r) $
in a suitable linear coordinate system),
where $ r : = m ( a+c-2) $.

Let $B''$ be the branch curve of a generic perturbation of $ p$:
then the braid monodromy factorization
corresponding to $B''$ can be calculated from the vertical and horizontal
braid
monodromies put together.

The problem of calculating the braid monodromy factorization corresponding
instead to the (cuspidal) branch curve $\De_m$ starting from the
braid monodromy factorization of $B''$ has been addressed, in greater
generality but in the  special case $ m=2$,
by Auroux and Katzarkov (\cite{a-k2}). Iteration of
their formulae should lead to the calculation of
   the braid monodromy factorization corresponding
to the (cuspidal) branch curve $\De_m$ in the case, sufficient for 
applications,
   where $m$ is a sufficiently large power of $2$.

Whether these formidable calculations will yield factorizations whose
m-equivalence is for us decidable  is still an
open question: but in both directions the result would be extremely
interesting, leading either to

i) a counterexample to the speculation DEF = CAN. SYMPL also in the simply
connected case, or to

ii) examples of diffeomorphic but not canonically symplectomorphic
simply connected algebraic surfaces.

\vfill

\noindent
{\bf Author's address:}

\bigskip

\noindent
Prof. Fabrizio Catanese\\
Lehrstuhl Mathematik VIII\\
Universit\"at Bayreuth, NWII\\
     D-95440 Bayreuth, Germany

e-mail: Fabrizio.Catanese@uni-bayreuth.de

\end{document}